\documentclass[12pt,twoside]{amsart}
\usepackage{amsmath}
\usepackage{amssymb}
\usepackage{amscd}
\usepackage{latexsym}
\usepackage{enumerate}

\begin{document}

\renewcommand{\baselinestretch}{1.3}
\renewcommand{\theequation}{\thesection.\arabic{equation}}
\newtheorem{theorem}{Theorem}[section]
\newtheorem{proposition}[theorem]{Proposition}
\newtheorem{corollary}[theorem]{Corollary}
\newtheorem{lemma}[theorem]{Lemma}
\newtheorem{mlemma}[theorem]{Main Lemma}
\newtheorem{mtheorem}[theorem]{Main Theorem}
\newtheorem{blackbox}[theorem]{Black Box}
\newtheorem{example}[theorem]{Example}
\newtheorem{remark}[theorem]{Remark}
\newtheorem{fact}[theorem]{Fact}
\newtheorem{construction}[theorem]{Construction}
\newtheorem{question}[theorem]{Question}
\newtheorem{definition}[theorem]{Definition}
\newtheorem{notation}[theorem]{Notation}
\newtheorem{titulo}[theorem]{}

\newenvironment{prooff}{\noindent{\sc Proof.}}{\qed\medskip}

\def\bull{{\unskip\nobreak\hfil\penalty50\hskip .001pt \hbox{}\nobreak\hfil
          \vrule height 1.2ex width 1.1ex depth -.1ex
           \parfillskip=0pt\finalhyphendemerits=0\medbreak}\rm}
\def\qed{\bull}
\def\End{\mathop{\rm End}\nolimits}
\def\dom{\mathop{\rm dom}\nolimits}
\def\Ext{\mathop{\rm Ext}\nolimits}
\def\Bext{\mathop{\rm Bext}\nolimits}
\def\Hom{\mathop{\rm Hom}\nolimits}
\def\Aut{\mathop{\rm Aut}\nolimits}
\def\ord{\mathop{\rm ord}\nolimits}
\def\sup{\mathop{\rm sup}\nolimits}
\def\dom{\mathop{\rm dom}\nolimits}

\def\Z{{\mathbb{Z}}}
\def\N{{\mathbb{N}}}
\def\P{{\mathbb{P}}}
\def\Q{{\mathbb{Q}}}

\def\ra{\longrightarrow}
\def\tT{\tilde{T}}
\def\tY{\tilde{Y}}
\def\tg{\tilde{g}}
\def\tt{\tilde{t}}
\def\tm{\tilde{m}}
\def\tz{\tilde{z}}
\def\tx{\tilde{x}}
\def\th{\tilde{h}}
\def\te{\tilde{\eta}}
\def\tpsi{\tilde{\psi}}
\def\tG{\tilde{G}}
\def\tH{\tilde{H}}
\def\tp{\tilde{\varphi}}


\title{It is consistent with ZFC that $B_1$-groups are not $B_2$}

\author[Saharon Shelah]{Saharon Shelah}
\address{Department of Mathematics, Hebrew University\\
91904 Jerusalem, Israel\\
and Rutgers University\\
Newbrunswick, NJ U.S.A.\\
e-Mail: Shelah@math.huji.ae.il\\}
\author[Lutz Str\"ungmann]{Lutz Str\"ungmann}
\address{Department of Mathematics, Hebrew University\\ 91904 Jerusalem,
Israel\\e-mail: lutz@math.huji.ac.il}
\thanks{
The first author was supported by project No. G-0545-173,06/97 of
the {\em German-Israeli Foundation for Scientific Research \&
Development}.\\ Publication 754 in Shelah's list of publications\\
The second author was supported by a MINERVA fellowship.}


\begin{abstract}
A torsion-free abelian group $B$ of arbitrary rank is called a
$B_1$-group if $\Bext^1(B,T)=0$ for every torsion abelian group
$T$, where $\Bext^1$ denotes the group of equivalence classes of
all balanced exact extensions of $T$ by $B$. It is a
long-standing problem whether or not the class of $B_1$-groups
coincides with the class of $B_2$-groups. A torsion-free abelian
group $B$ is called a $B_2$-group if there exists a continuous
well-ordered ascending chain of pure subgroups, $0=B_0 \subset B_1
\subset \cdots \subset B_{\alpha} \subset \cdots \subset
B_{\lambda}=B=\bigcup\limits_{\alpha \in \lambda}B_{\alpha}$ such
that $B_{\alpha +1}=B_{\alpha}+G_{\alpha}$ for every $\alpha \in
\lambda$ for some finite rank Butler group $G_{\alpha}.$ Both,
$B_1$-groups and $B_2$-groups are natural generalizations of
finite rank Butler groups to the infinite rank case and it is
known that every $B_2$-group is a $B_1$-group. Moreover, assuming
$V=L$ it was proven that the two classes coincide. Here we
demonstrate that it is undecidable in ZFC whether or not all
$B_1$-groups are $B_2$-groups. Using Cohen forcing we prove that
there is a model of ZFC in which there exists a $B_1$-group that
is not a $B_2$-group.

\end{abstract}

\maketitle


\section{Introduction}
The study of Butler groups, both in the finite and in the infinite
rank case, is a most active area of Abelian Group Theory. There
are several challenging problems which require deep insight into
the theory of Butler groups and the available methods as well as
the development of new machinery. The finite rank case is closely
related to the study of representations of finite posets while the
infinite rank case has its own special flavor. During the last
years more and more the connection between infinite rank Butler
groups and infinite combinatorics was discovered and led to
numerous interesting results. In this paper we discuss one of the
long-standing problems, namely whether or not all $B_1$-groups are
$B_2$-groups, and show that its solution is independent of ZFC.

All groups in the following are abelian. Recall that a
torsion-free group $B$ of finite rank is called a {\it Butler
group} if it is a pure subgroup of a completely decomposable
group of finite rank. Butler \cite{B} introduced this class of
groups and proved that being Butler is equivalent to being an
epimorphic image of a completely decomposable group of finite
rank. Later Bican and Salce \cite{BS} noticed that a torsion-free
group $B$ of finite rank is a Butler group if and only if
$\Bext^1(B,T)=0$ for all torsion groups $T$, where $\Bext^1$
denotes the group of equivalence classes of all balanced exact
extensions of $T$ by $B$. This result initiated a generalization
of Butler groups to the infinite rank case in two "different"
ways: a torsion-free group $B$ of arbitrary rank is called
\begin{enumerate}
\item a {\it $B_1$-group} if $\Bext^1(B,T)=0$ for all torsion groups $T$;
\item a {\it $B_2$-group} if there exists a continuous well-ordered ascending
chain of pure subgroups, \[ 0=B_0 \subset B_1 \subset \cdots
\subset B_{\alpha} \subset \cdots \subset
B_{\lambda}=B=\bigcup\limits_{\alpha \in \lambda}B_{\alpha} \]
such that $B_{\alpha +1}=B_{\alpha}+G_{\alpha}$ for every $\alpha
\in \lambda$ for some finite rank Butler group $G_{\alpha}.$
\end{enumerate}
The authors put the word different in quotation marks on purpose
since it was not known whether or not the two definitions 1. and
2. describe the same class of torsion-free groups. In fact it is
the main aim of this paper to show that it is undecidable in ZFC
whether or not $B_1$-groups and $B_2$-groups coincide. This
problem is just one among the major questions in the theory of
infinite rank Butler groups but we will not touch upon the others
in this paper.

Bican and Salce \cite{BS} proved that for countable groups the
two definitions of $B_1$-group and $B_2$-group coincide and
without any cardinality restriction, every $B_2$-group is a
$B_1$-group. That the two classes coincide even for groups up to
cardinality $\leq \aleph_1$ was observed by Dugas-Hill-Rangaswamy
\cite{DHR} and Albrecht-Hill \cite{AH}. However, for groups of
higher cardinality an affirmative answer needed additional
set-theory (e.g. assuming CH the two classes coincide for groups
up to cardinality $\leq \aleph_{\omega}$, see \cite{DHR}). Several
other results were obtained by using {\it $\aleph_0$-prebalanced}
chains and {\it Axiom-3 families} (see \cite{BF}, \cite{F2} and
\cite{FM}). A very nice result which uses an algebraic, rather
than a set-theoretic, condition to ensure that $B_1$-groups are
$B_2$-groups is the main result in Fuchs-Rangaswamy \cite{FR}:

{\it A $B_1$-group $B$ which is the union of a continuous
well-ordered ascending chain of pure subgroups $B_{\alpha}$ each
of which has a countable typeset is necessarily a $B_2$-group.}

Assuming the continuums hypothesis Rangaswamy \cite{R} showed
that a torsion-free group $B$ is a $B_2$-group if and only if
$\Bext^1(B,T)=\Bext^2(B,T)=0$ for all torsion groups $T$. Thus it
was natural to ask whether or not for every torsion-free group
$G$, $\Bext^2(G,T)=0$ for all torsion groups $T$. Under the
negation of CH a negative answer was given by Dugas-Thom\'{e}
\cite{DT} while Magidor and the first author \cite{MS} answered
this question to the negative even assuming the generalized
continuums hypothesis. One of the main results on Butler groups
of arbitrary cardinality assuming G\"odel's universe of
constructability (V=L) can be found in Fuchs-Magidor \cite{FM}:

{\it Assuming V=L every $B_1$-group is a $B_2$-group. }

In contrast to this result we will show that using Cohen forcing
there is a model of ZFC in which there exists a $B_1$-group that
is not a $B_2$-group. Hence it is undecidable in ZFC whether or
not the two classes of $B_1$-groups and $B_2$-groups coincide.

Our terminology is standard and maps are written on the left. If
$H$ is a subgroup of a torsion-free group then the purification
of $H$ in $G$ is denoted by $H_*$. For notations and basic facts
we refer to \cite{Fu} for abelian groups, \cite{K} and \cite{S}
for forcing and \cite{EM} or \cite{J} for set-theory. Moreover,
the interested reader may look at \cite{AV} for a survey on finite
rank Butler groups and at \cite{Bi}, \cite{F1} for surveys on
infinite rank Butler-groups.

\section{Infinite rank Butler groups}
In this section we recall the definitions of $B_1$-groups and
$B_2$-groups as they were given by Bican-Salce in \cite{BS}. Both
classes contain the class of finite rank Butler-groups (pure
subgroups of completely decomposable groups of finite rank) first
studied by Butler in \cite{B}. Let us begin with the notion of a
balanced subgroup.

A pure subgroup $A$ of the torsion-free group $G$ is said to be a
{\it balanced} subgroup if every coset $g+A$ ($g \in G$) contains
an element $g+a$ $(a \in A)$ such that $\chi(g+a) \geq \chi(g+x)$
for all $x \in A$. Such an element is called {\it proper with
respect to A} and $\chi(g)$ denotes the characteristic of an
element $g$ in the given group $G$.

An exact sequence $0 \rightarrow A \rightarrow G \rightarrow C
\rightarrow 0$ is {\it balanced exact} if the image of $A$ in $G$
is a balanced subgroup of $G$. Hunter \cite{H} discovered that
the equivalence classes of balanced extensions of a group $H$ by
a group $G$ give rise to a subfunctor $\Bext^1(H,G)$ of
$\Ext^1(H,G)$ and hence homological algebra is applicable. Thus
for a balanced exact sequence \begin{equation} 0 \rightarrow A
\rightarrow G \rightarrow C \rightarrow 0 \tag{$*$} \end{equation}
and a group $H$ we obtain the two long exact sequences
\[ 0 \rightarrow \Hom(C,H) \rightarrow \Hom(G,H) \rightarrow \Hom(A,H)
\rightarrow \Bext^1(C,H) \rightarrow \Bext^1(G,H) \rightarrow\] \[
\rightarrow \Bext^1(A,H) \rightarrow \Bext^2(C,H) \rightarrow
\cdots \] and
\[ 0 \rightarrow \Hom(H,A) \rightarrow \Hom(H,G) \rightarrow \Hom(H,C)
\rightarrow \Bext^1(H,A) \rightarrow \Bext^1(H,G) \rightarrow \]
\[ \rightarrow \Bext^1(H,C) \rightarrow \Bext^2(H,A) \rightarrow \cdots \]
It is routine to check that balanced-exactness of the sequence
$(*)$ is equivalent to the following property: for every rank $1$
torsion-free group $R$, every homomorphism $R \rightarrow C$ can
be lifted to a map $R \rightarrow G$, i.e. every rank $1$
torsion-free group is projective with respect to $(*)$. Thus the
following lemma is easily established.

\begin{lemma}
\label{local} Let
\[0 \rightarrow A \rightarrow G \overset{\varphi}{\rightarrow} C \rightarrow 0 \]
be a balanced exact sequence. Then this sequence is locally
invertible, i.e. for any element $c \in C$ there exists a
homomorphism $\psi_c : \left< c \right>_* \rightarrow G$ such that
$\varphi\psi_c=id_{\left< c \right>_*}$.
\end{lemma}

We now come to the definitions of $B_1$-groups and $B_2$-groups.

\begin{definition}
A torsion-free abelian group $B$ is called
\begin{enumerate}
\item a {\it $B_1$-group} if $\Bext^1(B,T)=0$ for all torsion groups $T$;
\item a {\it $B_2$-group} if there exists a continuous well-ordered ascending
chain of pure subgroups, \[ 0=B_0 \subset B_1 \subset \cdots
\subset B_{\alpha} \subset \cdots \subset
B_{\lambda}=B=\bigcup\limits_{\alpha \in \lambda}B_{\alpha} \]
such that $B_{\alpha +1}=B_{\alpha}+G_{\alpha}$ for every $\alpha
\in \lambda$ for some finite rank Butler group $G_{\alpha};$ i.e.
$B_{\alpha}$ is {\it descent} in $B_{\alpha+1}$ in the sense of
Albrecht-Hill \cite{AH};
\item {\it finitely Butler} if every finite rank subgroup of $B$ is
a Butler-group.
\end{enumerate}
\end{definition}

Due to Bican-Salce \cite{BS} the three definitions are equivalent
for countable torsion-free groups.

\begin{theorem}[\cite{BS}]
For a countable torsion-free abelian group $B$ the following are
equivalent:
\begin{enumerate}
\item $B$ is finitely Butler;
\item $B$ is a $B_2$-group;
\item $B$ is a $B_1$-group.
\end{enumerate}
\end{theorem}

Without any restriction to the cardinality we have in general:

\begin{theorem}[\cite{BS}]
$B_2$-groups of any rank are $B_1$-groups and finitely Butler.
\end{theorem}

It turned out that the converse implication in the above theorem
couldn't be proved without any additional set-theoretic
assumptions. There are some partial results in ZFC (mentioned in
the introduction) characterizing the $B_2$-groups among the
$B_1$-groups but non of them is really satisfactory. The
following was shown by Fuchs and Rangaswamy independently.

\begin{lemma}[\cite{F2}, \cite{R}]
Suppose that $0 \rightarrow H \rightarrow C \rightarrow G
\rightarrow 0$ is a balanced-exact sequence where $C$ is a
$B_2$-group and $H$ and $G$ are $B_1$-groups. If one of $H$ and
$G$ is a $B_2$-group, then so is the other. \end{lemma}

An attempt to characterize the $B_2$-groups in a homological way
is the following theorem due to Fuchs.

\begin{theorem}[\cite{F2}]
If $B$ is a $B_2$-group, then $\Bext^i(B,T)=0$ for all $i \geq 1$
and for all torsion groups $T$.
\end{theorem}

Assuming the continuums hypothesis Rangaswamy was able to show
that also the converse holds and in some cases Fuchs could even
remove CH.

\begin{theorem}[\cite{R}, \cite{F2}] The following is true:
\begin{enumerate}
\item Assuming CH a torsion-free group $B$ is a $B_2$-group if and
only if $\Bext^1(B,T)=\Bext^2(B,T)=0$ for all torsion groups $T$.
\item A torsion-free group $B$ of cardinality $\aleph_n$ (for an
integer $n \geq 1$) is a $B_2$-group if and only if
$\Bext^i(B,T)=0$ for all $i \leq n+1$ and all torsion groups $T$.
\end{enumerate}
\end{theorem}

Motivated by this result it was natural to ask whether
$\Bext^2(B,T)$ is always zero for a torsion-free group $B$ and a
torsion group $T$ but Magidor-Shelah \cite{MS} proved that this
is not the case even assuming the generalized continuums
hypothesis GCH. That CH was relevant in many papers was explained
by Fuchs who showed the following theorem.

\begin{theorem}[\cite{F2}]
In any model of ZFC, the following are equivalent:
\begin{enumerate}
\item $\Bext^2(G,T)=0$ for all torsion-free groups $G$ and torsion
groups $T$;
\item CH holds and balanced subgroups of completely decomposable
groups are $B_2$-groups.
\end{enumerate}
\end{theorem}

One of the most interesting and main results in the theory of
infinite rank Butler groups is the following final theorem of
this section proved by Magidor and Fuchs.

\begin{theorem}[\cite{FM}]
Assuming $V=L$ every $B_1$-group is a $B_2$-group.
\end{theorem}

We will show in this paper that the last theorem does not hold in
ZFC but is independent of ZFC.


\section{The forcing}

In this section we will explain the forcing notion we are going
to use to construct our $B_1$ group $H$ which fails to be $B_2$.
The reader who is familiar with forcing, especially with adding
Cohen reals may skip this section. Most results are well-known
and basic and for unexplained notations and for further results
on forcing we
refer to Kunen's book \cite{K} or more advanced to
Shelah's book \cite{S}.

Let $M$ be any countable transitive model of $ZFC$ and assume of
course that the set theory $ZFC$ is consistent. The aim of
forcing is to extend $M$ to a new model which still satisfies
$ZFC$ but which has additional properties which we are interested
in. Therefore we use the method of forcing which was first
discovered by Cohen.

A {\it forcing notion} $\P \in M$ is just a non empty, partially
ordered set $(\P, \leq, 0_{\P})$, where $0_{\P}$ is the minimal
element of $\P$, hence $0_{\P} \leq p$ for all $p \in \P$. Note
that we don't require that $p \leq q$ and $q \leq p$ implies
$q=p$. If two elements $p,q \in \P$ have no common upper bound,
i.e. there is no $t \in \P$ such that $q \leq t $ and $q \leq t$,
then we say that $p$ and $q$ are {\it incompatible} and write $p
\perp q$. If a common upper bound exists we call the elements
{\it compatible}. We now want to add to $M$ a subset $G$ of $\P$
to construct a transitive set $M[G]$ which is a model of $ZFC$
with the same ordinals as $M$ such that $M \subseteq M[G]$ and $G
\in M[G]$. Those sets $G$ are called generic.

\begin{definition}
Let $D \subseteq \P$, $G \subseteq \P$ and $p \in \P$. Then
\begin{enumerate}
\item $D$ is called {\rm dense} in $\P$ if for any $q \in \P$ there is
an element $t \in D$ such that $q \leq t$;
\item $D$ is {\rm dense above p} if for any $q \in \P$, $p \leq q$ there
exists an element $t \in D$ such that $q \leq t$;
\item $G$ is called {\rm $\P$-generic over $M$} if the following hold:
\begin{enumerate}
\item for all $q,r \in G$ there exists $t \in G$ such that $q \leq t$ and $r \leq t$,
i.e. all elements of $G$ are compatible;
\item if $q \in G$ and $t \leq q$ for some $t \in \P$ then also $t \in G$;
\item $G \cap D \not= \emptyset$ for every dense subset $D$ of $\P$ which is in $M$.
\end{enumerate}
\end{enumerate}
\end{definition}

A first observation is that a generic $G$ intersects also with
``dense above $p$'' sets in many cases.

\begin{lemma}
\label{denseabove}
Let $D \subseteq \P$ and $G$ be $\P$-generic
over $M$. Then
\begin{enumerate}
\item Either $G \cap D \not= \emptyset$ or there exists $q \in G$ such that for all $r \in D$
we have $r \perp q$;
\item If $p \in G$ and $D$ is dense above $p$, then $G \cap D \not= \emptyset$.
\end{enumerate}
\end{lemma}
\proof See \cite[Lemma 2.20]{K}. \qed

If $G$ is $\P$-generic over $M$ or for short generic, then the
existence of the model $M[G]$ with the desired properties follows
from the Forcing Theorem (see \cite{S}). $M[G]$ is the smallest
transitive model of $ZFC$ that contains $M$ and $G$. We don't
want to recall the construction of $M[G]$ but we would like to
mention the following facts. Since we want to prove theorems in
$M[G]$ we would like to know the members of $M[G]$ but we can not
have full knowledge of them inside $M$ since this would cause
these sets to be in $M$ already. If $G$ is in $M$ then $M[G]$
gives nothing new, so we have to assume that $G$ is not in $M$
and this is the case in general as the following lemma shows.

\begin{lemma}
Let $G$ be $\P$-generic over $M$. If $\P$ satisfies the following condition
\begin{equation}
\forall \text{ }  p \in \P \text{ }  \exists \text{ } q, \text{
}r \in \P \text{ such that } p \leq q, p \leq r \text{ and } q
\perp r
\end{equation}
then $G \not\in M$.
\end{lemma}

\proof See \cite[Lemma 2.4]{K}. \qed

Nevertheless, every element $p$ of $\P$ can be a member of a generic set.

\begin{lemma}
Let $p \in \P$. Then there is a $G$ which is $\P$-generic over $M$ such that $p \in G$.
\end{lemma}

\proof See \cite[Lemma 2.3]{K}. \qed

Although we don't know the generic set $G$ we assume that we have
some prescription for building the members of $M[G]$ out of $M$
and $G$. These prescriptions are called {\it $\P$-names}, usually
denoted by $\tau$, and their {\it interpretation in $M[G]$} is
$\tau[G]$. For the exact definition of $\P$-names and their
interpretation we refer again to Kunen's book \cite{K} but let us
mention that the Strengthened Forcing Theorem (see \cite{S})
shows that
\[ M[G] = \{ \tau[G] : \tau \in M \text{ and } \tau \text{ is a } \P \text{-name } \}. \]
If we are talking about the $\P$-name of a special object $H$
from $M[G]$ without specifying $G$ then we will write $\tilde{H}$
instead of $H$ to avoid confusion but if $H$ is already in $M$,
then we omit the snake. Any sentence of our forcing language uses
the $\P$-names to assert something about $M[G]$ but the truth or
falsity of a sentence $\psi$ in $M[G]$ depends on $G$ in general.
If $p \in \P$, then we write $p \Vdash \psi$ and say {\it p forces
$\psi$} to mean that for all $G$ which are $\P$-generic over $M$,
if $p \in G$, then $\psi$ is true in $M[G]$. If $0_{\P} \Vdash
\psi$ then we just write $\Vdash_{\P} \psi$ which means that for
any generic $G$ the sentence $\psi$ is true in $M[G]$ since
$0_{\P}$ is always contained in $G$. Hence the elements of $\P$
provide partial information about objects in $M[G]$ but not all
information and if $p \leq q$ then $q$ contains more information
than $p$. It is amazing but it may be decided in $M$ whether or
not $p \Vdash \psi$ and whenever a sentence $\psi$ is true in
$M[G]$ then there is $p \in G$ such that $p \Vdash \psi$.

We now turn to the forcing of adding Cohen reals. Therefore we
specify $\P$ and let $\kappa$ be an uncountable cardinal. We put

\[ \P =\{ p \mid p \text{ is a function from } \kappa \times \omega \text{ to } 2 \text{ with finite domain }\} \]
\[= \{ p \mid p : \kappa \times \omega \longrightarrow 2,\text{ } \dom(p) \text{ finite } \} \]

The partial ordering of $\P$ is given by set theoretic inclusion,
i.e. two functions $p$ and $q$ satisfy $p \leq q$ if and only if
$q$ extends $p$ as a function. This forcing is called {\it
``adding $\kappa$ Cohen reals''} and the elements of $\P$ can
obviously be regarded as functions from $\kappa$ to $^{<
\omega}2$ which we will do in the sequel.

The next lemma shows why the forcing is called adding $\kappa$
Cohen reals.

\begin{lemma}
$\Vdash_{\P}$ ``There are at least $\kappa$ reals''.
\end{lemma}

\proof See \cite[Chapter I, Lemma 3.3]{S}. \qed

We will give the $\kappa$ Cohen reals $\P$-names, say
$\tilde{\eta}_{\alpha}$ for $\alpha \in \kappa$ and state some
basic properties of the Cohen reals. Note that a real is a
function from $\omega$ to $2=\{0,1\}$.

\begin{lemma}
\label{reals}
The following hold for $\alpha, \beta \in \kappa$:
\begin{enumerate}
\item $\Vdash_{\P}$ ``There are infinitely many $n \in \N$ such that $\tilde{\eta}_{\alpha}(n)=\tilde{\eta}_{\beta}(n)=1$'';
\item $\Vdash_{\P}$ ``There are infinitely many $n \in \N$ such that $\tilde{\eta}_{\alpha}(n)=\tilde{\eta}_{\beta}(n)=0$'';
\item $\Vdash_{\P}$ ``There are infinitely many $n \in \N$ such that $\tilde{\eta}_{\alpha}(n)\not=\tilde{\eta}_{\beta}(n)$''.
\end{enumerate}
\end{lemma}

\proof The proof of this fact is standard using a densitiy
argument.
\qed

Moreover, we have three more important facts.

\begin{lemma}
The following hold for $\P$.
\begin{enumerate}
\item $\P$ satisfies the c.c.c. condition, i.e. $\P$ has no uncountable subset of pairwise
incompatible members;
\item $\P$ preserves cardinals and cofinalities, i.e. if $\lambda$ is a cardinal in $M$,
then $\lambda$ is also a cardinal in $M[G]$ with the same
cofinality;
\item $\Vdash_{\P}$ `` $2^{\aleph_0}\geq \lambda$''. In particular, if
$\lambda^{\aleph_0}=\lambda$ in $M$, then $\Vdash_{\P}$ ``
$2^{\aleph_0}= \lambda$''.
\end{enumerate}
\end{lemma}

\proof See \cite[Chapter I, Lemma 3.8]{S}, \cite[Chapter I,
Theorem 4.1]{S} and  \cite[Theorem 5.10]{K}. \qed

Finally we would like to remark that our notation is the
''Jerusalem style'' of forcing notation like in \cite{S} but
differs from the notation for example in \cite{K}. In our partial
order $p \leq q$ means that $q$ contains more information than
$p$ does and not vice versa.


\section{Our $B_1$ group $H$}

Let $M$ be a countable transitive model of $ZFC$ in which the
generalized continuum hypothesis holds, i.e.
$2^{\kappa}=\kappa^+$ for all infinite cardinals $\kappa$.
Moreover, let $\kappa \geq \aleph_4$ be regular and let $\P$ be
the forcing of adding $\kappa$ Cohen reals. As we have seen in the
last section, $\P$ preserves cardinals and cofinalities and
$2^{\aleph_0}=\kappa$ in $M[G]$ for every generic $G$. Let
$\tilde{\eta}_{\alpha}$ denote the Cohen reals for $\alpha \in
\kappa$ and let $M^*$ be a model of the Cohen forcing extending
$M$.

For the definition of our group $H$ we choose independent
elements
\[ \{ x_n : n \in \omega \} \text{ and } \{ y_{\alpha} : \alpha <
\kappa \} \] and fix a countable set of natural prime numbers
\[ \{ p_n \in \Pi : n \in \omega \} \]
such that $p_n < p_m$ for $n < m$. Here $\Pi$ denotes the set of
all primes.

\begin{definition}
\label{defi}
Let $W=\bigoplus\limits_{n \in \omega} \Q x_n \oplus
\bigoplus\limits_{\alpha < \kappa}\Q y_{\alpha}$ be the rational
vector space and let $F=\bigoplus\limits_{n \in \omega} \Z x_n
\oplus \bigoplus\limits_{\alpha < \kappa}\Z y_{\alpha}$ be the
free abelian group generated by the $x_n$'s and $y_{\alpha}$'s.
Inside $W$ we define
\[ \tH=\left< F, p_n^{-1}(y_{\alpha} - x_n) : \alpha < \kappa, n \in \omega,
\tilde{\eta}_{\alpha}(n)=1 \right> \subseteq W^+ \] as a subgroup
of the additive group $W^+$ of $W$.
\end{definition}

Letting $H$ being the interpretation of $\tH$ in $M^*$ we can now
state our Main Theorem.

\begin{mtheorem}
\label{main} In the model $M^*$ the group $H$ is a $B_1$-group but
not a $B_2$-group. Hence it is consistent with ordinary set theory
ZFC that $B_1$-groups need not be $B_2$-groups.
\end{mtheorem}

The proof of the Main Theorem \ref{main} will be divided into two
parts. The first part is to show that $H$ is a $B_1$-group which
will be done in this section. Section $5$ will then consist of
proving that $H$ is not $B_2$.

\begin{theorem}
\label{main1} In the model $M^*$ the group $H$ is a $B_1$-group.
\end{theorem}

The proof of Theorem \ref{main1} takes the rest of this section
and consists of several steps.

\proof {\it (of Theorem \ref{main1})} To prove that $H$ is a
$B_1$-group we have to show that $\Bext(H,T)=0$ for any torsion
group $T$. Hence let
\begin{equation}
0 \ra \tT \overset{id}{\ra} \tG \overset{\tp}{\ra} \tH \ra 0
\label{seq}
\end{equation}
be a balanced exact sequence with $\tT$ torsion. Thus there
exists $r^* \in \P$ such that
\[ r^* \Vdash \text{''} 0 \ra \tT \overset{id}{\ra} \tG \overset{\tp}{\ra} \tH \ra
0\text{ is balanced exact.''} \] We choose preimages $\tg_{\alpha}
\in G$ of $y_{\alpha}$ under $\tp$ such that
$\tp(\tg_{\alpha})=y_{\alpha}$ for all $\alpha < \kappa$.
Similarly let $\tx_n \in G$ be a preimage for $x_n$ under $\tp$
for $n \in \omega$. Moreover, let
\[ A_{\alpha}= \{ n \in \omega : \te_{\alpha}(n)=1 \} \]
for $\alpha < \kappa$.

It is our aim to show that the balanced exact sequence (\ref{seq})
is forced to split, hence it is enough to prove that the
homomorphism $\tp$ is invertible, i.e. we have to find $\tpsi :
\tH \ra \tG$ such that $\tp\tpsi=id_{\tH}$. Therefore it is
necessary to find preimages of the generators of $\tH$ in $\tG$
such that equations satisfied in $\tH$ also hold in $\tG$. We
need the following definition.

\begin{definition}
Let $\alpha < \kappa$ and $\tt\in \tT$ arbitrary. Then the set
$R_{\alpha,\tt}$ is defined as
\[ R_{\alpha, \tt}=\{ n \in A_{\alpha} : \tg_{\alpha} - \tt - \tx_n
\text{ is not divisible by } p_n \}. \]
\end{definition}

We will now use a purely group theoretic argument to show that if
for every $\alpha < \kappa$ there is a $\tt_{\alpha} \in \tT$
such that $R_{\alpha, \tt_{\alpha}}$ is finite implies that $\tp$
is invertible.

\begin{lemma}
\label{R-empty}
Let $\alpha < \kappa$ and let $\tt \in \tT$ such
that $R_{\alpha, \tt}$ is finite. Then there exists $\tt_{\alpha}
\in \tT$ such that $R_{\alpha, \tt_{\alpha}}=\emptyset$.
\end{lemma}

\proof Since $R_{\alpha,\tt}$ is finite we may assume without loss
of generality that $R_{\alpha, \tt}$ has minimal cardinality.
Assume that $R_{\alpha, \tt}$ is not empty and fix $n \in
R_{\alpha, \tt}$. By the primary decomposition theorem we
decompose $\tT$ as
\[ \tT = \tT_{p_n} \oplus \tT^{\prime} \]
where $\tT_{p_n}$ denotes the $p_n$-primary component of $\tT$.
Since $n \in A_{\alpha}$ it follows that $p_n$ divides
$y_{\alpha}-x_n$, hence there exists $\tz \in \tG$ such that
\[ \tp(\tz)=p_n^{-1}(y_{\alpha}-x_n). \]
Thus
\[ (\tg_{\alpha}-\tt - \tx_n) - p_n \tz \in \tT = \tT_{p_n} \oplus
\tT^{\prime} \] and therefore there exist $\tt_0 \in \tT_{p_n}$
and $\tt_1 \in \tT^{\prime}$ such that
\[ (\tg_{\alpha}-\tt - \tx_n) - p_n \tz = \tt_0 + \tt_1. \]
Since $\tT^{\prime}$ is $p_n$ divisible we can write
$\tt_1=p_n\tt_2$ for some $\tt_2 \in \tT^{\prime}$. Hence
\[ (\tg_{\alpha}-\tt - \tx_n) - p_n (\tz - \tt_2) = \tt_0. \]
We let $\tt^{\prime}=\tt + \tt_0$ and will show that $R_{\alpha,
\tt^{\prime}}$ has smaller cardinality than $R_{\alpha, \tt}$ - a
contradiction. By the choice of $\tt^{\prime}$ we have
\[ (\tg_{\alpha} - \tt^{\prime} - \tx_n)=
\tg_{\alpha}-\tt-\tt_0-\tx_n=p_n(\tz-\tt_2) \] and hence $n
\not\in R_{\alpha, \tt^{\prime}}$. But on the other side, if $m
\not\in R_{\alpha, \tt}$, then $p_m$ divides
$(\tg_{\alpha}-\tt-\tx_m)$ and thus $p_m$ divides
$(\tg_{\alpha}-(\tt^{\prime}-\tt_0)-\tx_m)$. Since $p_n \not= p_m$
it follows that $p_m$ divides $\tt_0$ and therefore $p_m$ divides
$(\tg_{\alpha}-\tt^{\prime}-\tx_m)$. Hence $m \not\in R_{\alpha,
\tt^{\prime}}$ showing that $R_{\alpha, \tt^{\prime}}$ is
strictly smaller than $R_{\alpha, \tt}$. This finishes the proof.
\qed

\begin{lemma}
\label{invertible} Assume that for every $\alpha < \kappa$ there
exists $\tt_{\alpha} \in \tT$ such that $R_{\alpha,
\tt_{\alpha}}$ is finite. Then $\tp$ is invertible and hence the
sequence (\ref{seq}) splits.
\end{lemma}

\proof By Lemma \ref{R-empty} we may assume without loss of
generality that for every $\alpha < \kappa$ the set $R_{\alpha,
\tt_{\alpha}}$ is empty. Thus for each $n \in A_{\alpha}$ we can
find $\tz_{\alpha, n} \in \tG$ such that
\[ p_n \tz_{\alpha, n} = \tg_{\alpha} - \tx_n - \tt_{\alpha}. \]
We now define a homomorphism $\tpsi : \tH \ra \tG$ as follows:
\begin{enumerate}
\item $\tpsi(x_n)=\tx_n$ $(n \in \omega)$;
\item $\tpsi(y_{\alpha})=\tg_{\alpha} - \tt_{\alpha}$ $(\alpha < \kappa)$;
\item $\tpsi(p_n^{-1}(y_{\alpha}-x_n))=\tz_{\alpha, n}$ $(\alpha < \kappa,$ $n \in A_{\alpha})$.
\end{enumerate}
It is now easy to check and therefore left to the reader that 1.,
2. and 3. induce a well-defined homomorphism $\tpsi : \tH \ra
\tG$ satisfying $\tp\tpsi=id_{\tH}$. \qed

{\it (Continuation of the proof of Theorem \ref{main1})} Up to
now we haven't used any forcing and we haven't worked in the
model $M^*$ but in $M$ itself. By the above Lemma
\ref{invertible} it remains to show that there exists for every
$\alpha < \kappa$ an element $\tt_{\alpha} \in \tT$ such that the
set $R_{\alpha, \tt_{\alpha}}$ is finite. To ensure this the
forcing comes into the game.

We define for $\alpha \not= \beta < \kappa$ the pure subgroup
$\tH_{\alpha, \beta}=\left< y_{\beta} - y_{\alpha} \right>_*$ of
$\tH$. Since the sequence (\ref{seq}) is balanced exact Lemma
\ref{local} shows that there exist homomorphisms
\[ \tpsi_{\alpha, \beta} : \tH_{\alpha, \beta} \ra \tG \text{ such that }
\tp\tpsi_{\alpha, \beta}=id_{\tH_{\alpha, \beta}}. \] Let
$\th_{\alpha, \beta} = \tpsi_{\alpha, \beta}(y_{\beta}-y_{\alpha})
\in \tG$, hence
\[ \tt_{\alpha, \beta}= \th_{\alpha, \beta} - (\tg_{\beta}-\tg_{\alpha}) \in \tT. \]
Since $\tT$ is a torsion group we can find $\tm_{\alpha, \beta}
\in \omega$ such that
\[ \ord(\tt_{\alpha, \beta})= \tm_{\alpha, \beta}. \]
We can now easily show

\begin{fact}
$r^* \Vdash \text{'' If } n > \tm_{\alpha, \beta}, \text{ then }
p_n \text{ divides } (\tg_{\beta} - \tg_{\alpha}) \text{ for } n
\in A_{\alpha} \cap A_{\beta} \text{ ''} $
\end{fact}

\proof If $n > \tm_{\alpha, \beta}$, then $p_n > \tm_{\alpha,
\beta}$ follows since the primes $p_m$ are increasing. Therefore
$\gcd(p_n, \tm_{\alpha, \beta})=1$ and thus $p_n$ divides
$(\th_{\alpha, \beta} - (\tg_{\beta}-\tg_{\alpha}))$. Moreover,
$\th_{\alpha, \beta}=\tpsi_{\alpha, \beta}(y_{\beta}- y_{\alpha})$
is divisible by $p_n$ since $n \in A_{\alpha} \cap A_{\beta}$.
Hence $p_n$ divides $(\tg_{\beta}-\tg_{\alpha})$. \qed

Now let $r^* \leq r_{\alpha, \beta} \in \P$ such that $r_{\alpha,
\beta}$ forces some value $m_{\alpha, \beta}$ to $\tm_{\alpha,
\beta}$, i.e.
\[ r_{\alpha, \beta} \Vdash \text{ '' } \tm_{\alpha, \beta} =
m_{\alpha, \beta} \text{ '' }. \] Without loss of generality we
assume that $\beta \in \dom(r_{\alpha, \beta})$ for all $\alpha,
\beta$. Since all elements of $\P$ are functions from $\kappa$ to
$2$ with finite domain, we may write for some $n_{\alpha, \beta}
\in \omega$
\[ \dom(r_{\alpha, \beta}) = \{ \gamma_{(\alpha, \beta, 0)}, \cdots,
\gamma_{(\alpha, \beta, n_{\alpha, \beta})} \} \subset \kappa, \]
where $\gamma_{(\alpha, \beta, i)} < \gamma_{(\alpha, \beta, j)}$
if $i < j \leq n_{\alpha, \beta}$. We would like to apply the
$\Delta$-Lemma to the functions $r_{\alpha, \beta}$ to obtain a
$\Delta$-system but unfortunately the functions $r_{\alpha,
\beta}$ depend on two variables. This forces us to do the
$\Delta$-Lemma 'by hand'. For this we use the Erd\"os-Rado
Theorem (see \cite{ER}).

First we define a coloring on four tuples in $\aleph_1$. Let
$\alpha_0, \alpha_1, \alpha_2, \alpha_3 \in \aleph_1$ such that
$\alpha_0\leq \alpha_1 \leq \alpha_2 \leq \alpha_3$ and let
\[ c(\alpha_0, \alpha_1, \alpha_2, \alpha_3)\] consist of the
following entries in an arbitrary but fixed order:
\renewcommand{\labelenumi}{{\rm (}\roman{enumi}{\hspace*{0.1cm}\rm)}}
\begin{enumerate}
\item $n_{\alpha_0, \alpha_1}$;
\item $m_{\alpha_0, \alpha_1}$;
\item $r_{\alpha_0, \alpha_1}(\gamma_{(\alpha_0,\alpha_1,j)} : j \leq n_{\alpha_0,
\alpha_1})$;
\item $(tv(\gamma_{(\alpha_{n_1}, \alpha_{n_2}, n_3)} < \gamma_{(\alpha_{m_1},
\alpha_{m_2}, m_3)}) : n_1,n_2,m_1,m_2 < 4; n_3 <
n_{\alpha_{n_1}, \alpha_{n_2}}; m_3 < n_{\alpha_{m_1},
\alpha_{m_2}}).$
\end{enumerate}
Recall that $tv$ denotes the {\it truth-value} of the inequality
and we may assume without loss of generality that it is just Yes
or No. The above coloring is a coloring with less than $\omega$
colors and thus we may apply the Erd\"os-Rado Theorem. Note that
we are still working in our model $M$ in which $GCH$ holds by
assumption. Hence we have \[\aleph_4 \ra (\aleph_1)_{\aleph_0}^4\]
which is exactly what we need to apply the Erd\"os-Rado Theorem.
We obtain an increasing chain of $c$-homogeneous elements
\[ \Gamma = \{ \alpha_{\epsilon} : \epsilon < \omega_1+1 \} \]
which means that whenever $\alpha_{\epsilon_1},
\alpha_{\epsilon_2}, \alpha_{\epsilon_3}, \alpha_{\epsilon_4} \in
\Gamma$ such that $\alpha_{\epsilon_1}\leq \alpha_{\epsilon_2}\leq
\alpha_{\epsilon_3}\leq \alpha_{\epsilon_4}$ , then \[
c(\alpha_{\epsilon_1}, \alpha_{\epsilon_2}, \alpha_{\epsilon_3},
\alpha_{\epsilon_4})=c^* \] for a fixed color $c^*$. Let this
particular color consist of the following entries:

\renewcommand{\labelenumi}{{\rm (}\Roman{enumi}{\hspace*{0.1cm}\rm)}}
\begin{enumerate}
\item $n^*$;
\item $m^*$;
\item $(k_1,\cdots, k_{n^*})$ $(k_i \in \{0,1\})$;
\item $(l_1,\cdots, l_{16^2(n^*)^2})$ $(l_i \in \{$Yes, No$\})$.
\end{enumerate}

Let us first explain what the homogenity implies. Let
$\alpha_{\epsilon_1}, \alpha_{\epsilon_2} \in \Gamma$ such that
$\epsilon_1< \epsilon_2$, then (I) ensures that the domain of
$r_{\alpha_{\epsilon_1}, \alpha_{\epsilon_2}}$ has size $n^*$.
Moreover, (II) says that $r_{\alpha_{\epsilon_1},
\alpha_{\epsilon_2}}$ forces the value $m^*$ to
$m_{\alpha_{\epsilon_1}, \alpha_{\epsilon_2}}$ and (III) implies
that the image of $r_{\alpha_{\epsilon_1}, \alpha_{\epsilon_2}}$
is uniquely determined. Finally (IV) ensures that if we take
another pair $\alpha_{\epsilon_3}, \alpha_{\epsilon_4} \in
\Gamma$ such that $\epsilon_3 < \epsilon_4$, then the relationship
between the elements of the domains of $r_{\alpha_{\epsilon_1},
\alpha_{\epsilon_2}}$ and $r_{\alpha_{\epsilon_3},
\alpha_{\epsilon_4}}$ is fixed.

In the sequel we need to be above all the "trouble", hence we may
assume without loss of generality that $m^*$ is greater or equal
to length$(r_{\alpha_{\epsilon},
\alpha_{\rho}}(\gamma_{\alpha_{\epsilon}, \alpha_{\rho},e}))$ for
all $\epsilon < \rho < \omega_1$ and $e \leq n^*$. We can now
approach to the $\Delta$-Lemma.

\begin{definition}
\label{defidelta}
For $\alpha_{\epsilon} \in \Gamma$ we define
\renewcommand{\labelenumi}{{\rm (}\arabic{enumi}{\hspace*{0.1cm}\rm)}}
\begin{enumerate}
\item $u_{\alpha_{\epsilon}} = \dom (r_{\alpha_{\epsilon},
\alpha_{\epsilon+1}}) \cap \dom (r_{\alpha_{\epsilon},
\alpha_{\epsilon+2}})$;
\item $u^*=\bigcap\limits_{\epsilon< \omega_1}
u_{\alpha_{\epsilon}}$;
\item $s_{\epsilon} = r_{\alpha_{\epsilon},
\alpha_{\epsilon+1}} \restriction_{u_{\alpha_{\epsilon}}} =
r_{\alpha_{\epsilon}, \alpha_{\epsilon+2}}
\restriction_{u_{\alpha_{\epsilon}}}$.
\end{enumerate}
\end{definition}

The reader might ask why (3) in the above Definition
\ref{defidelta} is well-defined. This follows from homogenity
since (II) implies that for $\gamma \in u_{\alpha_{\epsilon}}$ we
have $r_{\alpha_{\epsilon},
\alpha_{\epsilon+1}}(\gamma)=r_{\alpha_{\epsilon},
\alpha_{\epsilon+2}}(\gamma)$. We are now ready to show the
following lemma, our version of the $\Delta$-system. Note that if
we talk about a $\Delta$-system of functions then we mean that
the corresponding domains of the functions form a $\Delta$-system.

\begin{lemma}
\label{deltasystem}
For $\alpha_{\epsilon}, \alpha_{\rho} \in
\Gamma$ such that $\epsilon < \rho$ we have
\[ u_{\alpha_{\epsilon}} \cap u_{\alpha_{\rho}} = u^*. \]
Hence the functions $s_{\epsilon}$ ($\alpha_{\epsilon} \in
\Gamma$) form a $\Delta$-system with root $u^*$. Moreover, for
fixed $\epsilon < \omega_1$ the functions $r_{\alpha_{\epsilon},
\alpha_{\rho}}$ $(\epsilon < \rho < \omega_1)$ form a
$\Delta$-system with root $u_{\alpha_{\epsilon}}$.
\end{lemma}

\proof Let $\alpha_{\epsilon}, \alpha_{\rho} \in \Gamma$ such
that $\epsilon < \rho$. Clearly we have $u^* \subseteq
u_{\alpha_{\epsilon}} \cap u_{\alpha_{\rho}}$ by Definition
\ref{defidelta}. It remains to show the converse inclusion.
Therefore let $\gamma \in u_{\alpha_{\epsilon}} \cap
u_{\alpha_{\rho}}$ and choose $\tau < \omega_1$ arbitrary. We
have to prove that $\gamma$ lies in $u_{\alpha_{\tau}}$.

If $\tau=\epsilon$ or $\tau=\rho$, then we are done.

If $\tau \geq \epsilon+1$, then $c(\alpha_{\epsilon},
\alpha_{\epsilon+1}, \alpha_{\tau}, \alpha_{\tau+1})=c^*$ by
homogenity. Since $\gamma \in \dom(r_{\alpha_{\epsilon},
\alpha_{\epsilon+1}})$ we can find $i \leq n^*$ such that $\gamma
= \gamma_{(\alpha_{\epsilon}, \alpha_{\epsilon+1}, i)}$ and
similarly $\gamma = \gamma_{(\alpha_{\rho}, \alpha_{\rho+1}, j)}$
for some $j\leq n^*$. It follows now that
\[ tv(\gamma_{(\alpha_{\epsilon}, \alpha_{\epsilon+1}, i)} < \gamma_{(\alpha_{\rho}, \alpha_{\rho+1},
j)})= \text{ No and } tv(\gamma_{(\alpha_{\rho}, \alpha_{\rho+1},
j)} < \gamma_{(\alpha_{\epsilon}, \alpha_{\epsilon+1}, i)})=
\text{ No.}\]  Hence there exists by homogenity $k \leq n^*$ such
that
\[ tv(\gamma_{(\alpha_{\epsilon}, \alpha_{\epsilon+1}, i)} < \gamma_{(\alpha_{\tau}, \alpha_{\tau+1},
k)})= \text{ No and } tv(\gamma_{(\alpha_{\tau}, \alpha_{\tau+1},
k)} < \gamma_{(\alpha_{\epsilon}, \alpha_{\epsilon+1}, i)})=
\text{ No.}\] Thus $\gamma = \gamma_{(\alpha_{\tau},
\alpha_{\tau+1}, k)} \in \dom (r_{\alpha_{\tau},
\alpha_{\tau+1}})$. Similarly it follows that $\gamma \in \dom
(r_{\alpha_{\tau}, \alpha_{\tau+2}})$ and hence $\gamma \in
u_{\alpha_{\tau}}$ which was to prove.

If $\tau < \epsilon+1$, then we use similar arguments to those
above to prove that $\gamma \in u_{\alpha_{\tau}}$.

Thus we have shown that $\gamma \in u_{\alpha_{\tau}}$ for any
$\tau < \omega_1$ and therefore $\gamma \in u^*$.

The same kind of arguments show that also the functions
$r_{\alpha_{\epsilon}, \alpha_{\rho}}$ $(\epsilon < \rho <
\omega_1)$ form a $\Delta$-system with root
$u_{\alpha_{\epsilon}}$ for fixed $\epsilon < \omega_1$. \qed

It is now easy to see by a pigeon-hole argument that we may
assume without loss of generality (and we will assume this in the
sequel) that all the functions from the $\Delta$-systems in Lemma
\ref{deltasystem} coincide on their root.

{\it (Continuation of the proof of Theorem \ref{main1})} The
following definition now makes sense.

\begin{definition}
For $\epsilon < \rho < \omega_1$ and a generic $G \subseteq \P$
we define
\begin{enumerate}
\item $s^*=s_{\epsilon}
\restriction_{u^*}=s_{\epsilon}\restriction_{(u_{\alpha_{\epsilon}}
\cap u_{\alpha_{\rho}})}$;
\item $\tY = \{ \alpha_{\tau} : s_{\tau} \in G \}$.
\end{enumerate}
\end{definition}

We can now show that $s^*$ is strong enough to force that $\tY$
has cardinality $\aleph_1$.

\begin{fact}
$s^* \Vdash \text{ ''} \mid \tY \mid = \aleph_1 \text{ ''} $
\end{fact}

\proof Let $\tY^{\prime}=\{ s_{\epsilon} : \alpha_{\epsilon} \in
\Gamma \backslash \tY \}$ and assume that $s^*$ does not force
$\tY$ to be of size $\aleph_1$. Then $\mid \tY^{\prime} \mid =
\aleph_1$. We will show that this set is pre-dense above $s^*$.
Therefore let $f \geq s^*$, then $\dom(f)$ is a finite subset of
$\kappa$ and $u^* \subseteq \dom(f)$. We choose $s_{\epsilon} \in
\tY^{\prime}$ such that $\dom(s_{\epsilon}) \backslash u^*$ is
disjoint to $\dom(f) \backslash u^*$. This is possible since by
Lemma \ref{deltasystem} the $s_{\epsilon}$'s form a
$\Delta$-system, hence
\[ \dom(s_{\tau}) \backslash u^* \cap \dom(s_{\beta}) \backslash
u^* = \emptyset \] for $\beta \not= \tau$. Now, $f$ and
$s_{\epsilon}$ are compatible and thus $\tY^{\prime}$ is pre-dense
above $s^*$. Therefore $\tY^{\prime} \cap G \not= \emptyset$ by
Lemma \ref{denseabove} - a contradiction. \qed

We are almost done now and prove the following statement.

\begin{fact}
\label{lastbutfinal} $s^* \Vdash \text{ '' If } \alpha_{\epsilon},
\alpha_{\rho} \in \tY \text{ and } n \in A_{\alpha_{\epsilon}}
\cap A_{\alpha_{\rho}} \backslash [0,m^*] \text{ then } p_n \text{
divides } g_{\alpha_{\epsilon}} - g_{\alpha_{\rho}} \text{ ''}$.
\end{fact}

\proof Let $s^* \leq s$ such that
\[ s \Vdash \text{ ''} n \in A_{\alpha_{\epsilon}} \cap
A_{\alpha_{\rho}} \backslash [0,m^*] \text{ ''}. \] Without loss
of generality we may assume that $s$ also forces truth values to
$\alpha_{\epsilon} \in \tY$ and $\alpha_{\rho} \in \tY$. If one
of them is No, then we are done and hence let us assume that both
are Yes. We will show that there exists $\gamma < \omega_1$ such
that
\begin{enumerate}
\item $\gamma > \epsilon$;
\item $\gamma > \rho$;
\item $\dom(r_{\alpha_{\epsilon}, \alpha_{\gamma}}) \backslash
u_{\alpha_{\epsilon}} \cup \dom(r_{\alpha_{\rho},
\alpha_{\gamma}})\backslash u_{\alpha_{\rho}} \cup \{
\alpha_{\gamma} \} \cup u_{\alpha_{\gamma}} \backslash u^*$ is
disjoint to $\dom(s)$.
\end{enumerate}
Obviously we can choose $\gamma > \epsilon, \rho$ such that
$\dom(s)$ is disjoint to $\{ \alpha_{\gamma} \}$, so all we have
to ensure is that also $\dom(r_{\alpha_{\epsilon},
\alpha_{\gamma}}) \backslash u_{\alpha_{\epsilon}} \cup
\dom(r_{\alpha_{\rho}, \alpha_{\gamma}})\backslash
u_{\alpha_{\rho}} \cup u_{\alpha_{\gamma}} \backslash u^*$ is
disjoint to $\dom(s)$. For this we prove that the three sets
\renewcommand{\labelenumi}{{\rm (}\arabic{enumi}{\hspace*{0.1cm}\rm)}}
\begin{enumerate}
\item $\{ \gamma < \omega_1 : \dom(r_{\alpha_{\epsilon},
\alpha_{\gamma}}) \backslash u_{\alpha_{\epsilon}}$ is not
disjoint to $\dom(s) \}$;
\item $\{ \gamma < \omega_1 : \dom(r_{\alpha_{\rho}, \alpha_{\gamma}})\backslash
u_{\alpha_{\rho}}$ is not disjoint to $\dom(s) \}$;
\item $\{ \gamma < \omega_1 : u_{\alpha_{\gamma}} \backslash u^*$
is not disjoint to $\dom(s) \}$.
\end{enumerate}
are bounded in $\omega_1$. Let us start with $(1)$. By Lemma
\ref{deltasystem} we know that for each $\epsilon < \omega_1$ the
domains $\{ \dom(r_{\alpha_{\epsilon}, \alpha_{\gamma}}) :
\epsilon < \gamma < \omega_1 \}$ form a $\Delta$-system with root
$u_{\alpha_{\epsilon}}$, hence $\{ \dom(r_{\alpha_{\epsilon},
\alpha_{\gamma}})\backslash u_{\alpha_{\epsilon}} : \epsilon <
\gamma < \omega_1 \}$ is a set of pairwise disjoint sets. Since
$\dom(s)$ is a finite set $\{ \gamma < \omega_1 :
\dom(r_{\alpha_{\epsilon}, \alpha_{\gamma}}) \backslash
u_{\alpha_{\epsilon}}$ is not disjoint to $\dom(s) \}$ must be
bounded in $\omega_1$. Similarly $\{ \gamma < \omega_1 :
\dom(r_{\alpha_{\rho}, \alpha_{\gamma}})\backslash
u_{\alpha_{\rho}}$ is not disjoint to $\dom(s) \}$ is bounded in
$\omega_1$. Finally, again by Lemma \ref{deltasystem} the sets
$\{ u_{\alpha_{\gamma}} : \gamma < \omega_1 \}$ form a
$\Delta$-system with root $u^*$ and so also $\{ \gamma < \omega_1
: u_{\alpha_{\gamma}} \backslash u^*$ is not disjoint to $\dom(s)
\}$ is bounded in $\omega_1$.

For this $\gamma$ we are able to prove that there is $s^+ \geq s$
such that
\renewcommand{\labelenumi}{{\rm (}\roman{enumi}{\hspace*{0.1cm}\rm)}}
\begin{enumerate}
\item $s^+ \geq s$;
\item $s^+ \Vdash \text{''} \eta_{\alpha_{\gamma}}(n)=1 \text{''}$;
\item $s^+ \geq r_{\alpha_{\epsilon, \gamma}}$;
\item $s^+ \geq r_{\alpha_{\rho, \gamma}}$.
\end{enumerate}

Since $n$ was chosen large enough which means that
$\eta_{\alpha_{\gamma}}$ has length less or equal to $m^*$ and
hence less or equal to $n$, there is, once we know that we can
satisfy $(i)$, $(iii)$ and $(iv)$, also some $s^+ \geq s$
satisfying all conditions $(i)$, $(ii)$, $(iii)$ and $(iv)$. Thus
we only have to satisfy conditions $(i)$, $(iii)$ and $(iv)$ and
for this it is obviously enough to show that the three functions
$s$, $r_{\alpha_{\epsilon},\alpha_{\gamma}}$ and
$r_{\alpha_{\rho},\alpha_{\gamma}}$ are compatible. Assume that
$r_{\alpha_{\epsilon},\alpha_{\gamma}}$ and
$r_{\alpha_{\rho},\alpha_{\gamma}}$ are incompatible, then by
induction we obtain that
$r_{\alpha_{\epsilon},\alpha_{\omega_1}}$ and
$r_{\alpha_{\rho},\alpha_{\omega_1}}$ are incompatible. Hence for
$\tau < \sigma < \omega_1$ we have that
$r_{\alpha_{\tau},\alpha_{\omega_1}}$ and
$r_{\alpha_{\sigma},\alpha_{\omega_1}}$ are incompatible which
contradicts the c.c.c. condition of our forcing. Therefore
$r_{\alpha_{\epsilon},\alpha_{\gamma}}$ and
$r_{\alpha_{\rho},\alpha_{\gamma}}$ are compatible. Finally $s$
and $r_{\alpha_{\epsilon},\alpha_{\gamma}}$ (and similarly $s$
and $r_{\alpha_{\rho},\alpha_{\gamma}}$) are compatible since by
the choice of $\gamma$ we have $\dom(s) \cap
\dom(r_{\alpha_{\epsilon},\alpha_{\gamma}}) =
u_{\alpha_{\epsilon}}$.

Now $s^+ \Vdash \text{''} p_n$ divides $g_{\alpha_{\epsilon}} -
g_{\alpha_{\gamma}} \text{''}$ and $s^+ \Vdash \text{''} p_n$
divides $g_{\alpha_{\gamma}} - g_{\alpha_{\rho}} \text{''}$ and
therefore $s^+ \Vdash \text{''} p_n$ divides
$g_{\alpha_{\epsilon}} - g_{\alpha_{\rho}} \text{''}$ as claimed.
\qed

Finally we have to prove another fact.

\begin{fact}
\label{finalfact} $s^* \Vdash \text{ '' For every } \beta < \kappa
\text{ there exists } m_{\beta} < \omega \text{ such that }
\forall m_{\beta}< n \in A_{\beta} \text{ and } \alpha_{\epsilon}
\in \tY \text{ such that } n \in A_{\alpha_{\epsilon}} \text{ we
have } p_n \text{ divides } g_{\beta} - g_{\alpha_{\epsilon}}
\text{ ''}. $
\end{fact}

First note that this implies that the set $R_{\beta, 0}$ is
contained in $[0,m_{\beta})$ for all $\beta < \kappa$ and hence
finite after modifying the choice of the preimages $\tx_n$ of
$x_n$ $(n \in \omega)$ slightly (which doesn't has any effect on
what we have done so far). Choose $\tx_n \in \tG$ such that
\begin{enumerate}
\item $\varphi(\tx_n)=x_n$;
\item if $n > m^*$ and $\alpha \in \tY$ such that
$\eta_{\alpha}(n)=1$, then $p_n$ divides $g_{\alpha} - \tx_n$.
\end{enumerate}

For example if $n > m^*$ we choose $\alpha \in \tY$ such that
$\eta_{\alpha}(n)=1$; Let $k_n \in \tG$ such that $\varphi(k_n)=
1/p_n(y_{\alpha}-x_n)$ and put $\tx_n=p_nk_n + g_{\alpha}$. Then
clearly (i) and (ii) are satisfied.

Now Fact \ref{finalfact} ensures that $R_{\beta,0}$ is contained
in $[0, m_{\beta})$ because: For $m_{\beta} < n \in A_{\beta}$
choose $\alpha_{\epsilon} \in \tY$ such that
$\eta_{\alpha_\epsilon}(n)=1$ (the one which was used when
choosing the $\tx_n$'s), then we have by the choice of $\tx_n$
that $p_n$ divides $g_{\alpha_\epsilon} - \tx_n$ and by Fact
\ref{finalfact} we have $p_n$ divides
$g_{\beta}-g_{\alpha_\epsilon}$ and hence $p_n$ divides
$g_{\beta} - \tx_n$. Thus $n \not\in R_{[\beta,o)}$ and
$R_{[\beta,0)} \subseteq [0, m_{\beta})$ follows.

Therefore the proof of Fact \ref{finalfact} finishes the proof of
Theorem \ref{main1}.

\proof Fix $\beta < \kappa$ and let $s^+<s$ such that $s$ forces
$n \in A_{\beta}$. For every $\epsilon < \omega_1$ we choose (if
possible) $t_{\beta}^{\epsilon}$ in the generic set such that
\begin{enumerate}
\item $s \leq t_{\beta}^{\epsilon}$;
\item $t_{\beta}^{\epsilon} \Vdash \text{''} \alpha_{\epsilon} \in
\tY \text{''}$;
\item $t_{\beta}^{\epsilon} \Vdash \text{''}
\tm_{\alpha_{\epsilon},\beta} =
m^{\sharp}_{\alpha_{\epsilon},\beta} \text{''}$ for some
$m^{\sharp}_{\alpha_{\epsilon}, \beta} \in \N$.
\end{enumerate}
Note that it is sufficient to find one $\epsilon$ s.t.
\begin{equation}
p_n \text{ divides } g_\beta - g_{\alpha_\epsilon} \tag{$*$}
\end{equation}
, for then we can use fact \ref{lastbutfinal} to get the
conclusion for any $\alpha_\rho$ s.t. $n \in A_{\alpha_\rho}$. If
we have one $t^\epsilon_\beta$ satisfying (ii) and (iii), then it
forces ($*$) for $n > m^{\sharp}_{\alpha_\epsilon , \beta}$. So we
first ensure (ii) and (iii) and then we use that there is an
uncountable subset $S_{\beta}$ of $\omega_1$ such that $\{
t_{\beta}^{\epsilon} : \epsilon \in S_{\beta} \}$ is a
$\Delta$-system to ensure (i) where we put
$m_{\beta}=m^{\sharp}_{\alpha_{\epsilon, \beta}}$ which can be
chosen fixed for the $\Delta$-system.
\qed


\section{Why $H$ fails to be $B_2$}
To complete the proof of our Main Theorem \ref{main} we show in
this section that the interpretation $H$ of the group $\tH$ from
Definition \ref{defi} can not be a $B_2$-group in $M^*$.

\begin{theorem}
\label{main2} In the model $M^*$ the group $H$ can not be a
$B_2$-group. \end{theorem}

\proof Towards contradiction assume that $H$ is a $B_2$-group,
hence has a $B_2$-filtration
\[ H=\bigcup\limits_{\alpha < \kappa} H_{\alpha}. \]
Recall that a $B_2$-filtration is a smooth ascending chain of
pure subgroups $H_{\alpha}$ such that for every $\alpha < \kappa$
$H_{\alpha + 1} = H_{\alpha} + B_{\alpha}$ for some finite rank
Butler group $B_{\alpha}$. We need the following lemma and recall
that a {\it cub} in $\kappa$ is a subset $C$ of $\kappa$ such that
\begin{enumerate}
\item $C$ is {\it closed} in $\kappa$, i.e. for all $C^{\prime}
\subseteq C$, if $\sup C^{\prime} < \kappa$, then $\sup
C^{\prime} \in C$;
\item $C$ is {\it unbounded} in $\kappa$, i.e. $\sup C^{\prime} = \kappa$.
\end{enumerate}

The proof of the following lemma is standard (see
\cite{EM}[II.4.12]) but for the convenience of the reader we
include it briefly.

\begin{lemma}
\label{cub}
The set $C=\{ \delta < \kappa \mid H_{\delta}=\left<
x_n, y_{\beta} : n \in \omega, \beta < \delta \right>_* \}$ is a
closed unbounded set (cub) in $\kappa$.
\end{lemma}

\proof First we show that $C$ is closed in $\kappa$. Therefore
let $C^{\prime}=\{ \delta_i \mid i \in I \}$ be a subset of $C$
such that $\sup C^{\prime} < \kappa$. If we put $\gamma =\sup
C^{\prime}$, then clearly
\[ H_{\gamma}=\bigcup\limits_{i \in I}H_{\delta_i}=\left< x_n,
y_{\beta} : n \in \omega, \beta < \gamma \right>_* \] and hence
$\gamma \in C$.

It remains to show that $C$ is unbounded. Therefore assume that
$C$ is bounded by $\delta^* < \kappa$, i.e. $\delta \leq \delta^*$
for all $\delta \in C$. We will show that there exists $\delta^* <
\gamma <\kappa$ such that $H_{\gamma}=\left< x_n, y_{\beta} : n
\in \omega, \beta < \gamma \right>_*$, hence $\gamma \in C$ - a
contradiction.

Let $\rho_1 = \delta^*$ and put
\[ E_{\rho_1} = \left< x_n, y_{\beta} : n \in \omega, \beta <
\rho_1 \right>_*.\] Now choose $\rho_1 \leq \alpha_1 < \kappa$
such that $E_{\rho_1} \subseteq H_{\alpha_1}$. If $\alpha_1
\not\in C$ choose $\alpha_1 \leq \rho_2 < \kappa$ such that
\[ H_{\alpha_1} \subseteq E_{\rho_2}=\left< x_n, y_{\beta} : n
\in \omega, \beta < \rho_2 \right>_*. \] Continuing this way we
obtain a sequence of groups $E_{\rho_i}$ and $H_{\alpha_i}$ such
that
\[ E_{\rho_i} \subseteq H_{\alpha_i} \subseteq E_{\rho_{i+1}}
\]
for all $i \in \omega$. Let $\gamma = \sup \{ \rho_i : i \in
\omega \} = \sup \{ \alpha_i : i \in \omega \}$. Then
\[ H_{\rho} = \bigcup\limits_{i \in \omega} H_{\alpha_i} =
\bigcup\limits_{i \in \omega} E_{\rho_i} = E_{\gamma}= \left<
x_n, y_{\beta} : n \in \omega, \beta < \gamma \right>_* \] and
hence $\gamma \in C$. This finishes the proof. \qed

{\it (Continuation of the proof of Theorem \ref{main2})} Now let
$\delta \in C$ such that $\delta > \aleph_1$. This is possible
since $C$ is a cub by the previous Lemma \ref{cub}. Note that
$y_{\delta} \not\in H_{\delta}$ but we have the following lemma.

\begin{lemma}
\label{contained} There exists $n^* \in \omega$ and a sequence of
ordinals $\delta \leq \alpha_1 \leq \alpha_2 \cdots \leq
\alpha_{n^*} < \kappa$ such that
\[ \left< H_{\delta} + \Z y_{\delta} \right>_* \subseteq
\sum\limits_{m \leq n^*}B_{\alpha_m} + H_{\delta}. \]
\end{lemma}

\proof We induct on $\alpha \geq \delta$ to show the even stronger
statement that for any $L \subseteq_* H_{\alpha}$, $L$ of finite
rank, there exist $n^* \in \omega$ and $\delta \leq \alpha_1 \leq
\alpha_2 \cdots \leq \alpha_{n^*} < \kappa$ such that
\[ \left< H_{\delta} + L \right>_* \subseteq
\sum\limits_{m \leq n^*}B_{\alpha_m} + H_{\delta}. \]

If $\alpha=\delta$, then we are done choosing $n^*=1$ and
$\alpha_1=\alpha$.

If $\alpha > \delta$ is a limit ordinal, then $L \subseteq_*
H_{\alpha}$ implies $L \subseteq_* H_{\beta}$ for some $\delta
\leq \beta < \alpha$. Hence we are done by induction hypothesis.

If $\alpha =\beta +1$, let $H_{\alpha} = H_{\beta} + B_{\beta}$
and let $L=\left< l_1, \cdots, l_k \right>_*$. We can find
representations
\[ l_i = h_{\beta,i} + b_{\beta, i} \]
for all $1 \leq i \leq k$ where $h_{\beta, i} \in H_{\beta}$ and
$b_{\beta, i} \in B_{\beta}$. We put
\[ L_{\beta}=\left< h_{\beta, i} , (B_{\beta} \cap H_{\beta}) : 1 \leq i \leq
k\right>_* \subseteq H_{\beta} \] which is a pure subgroup of
finite rank of $H_{\beta}$. An easy calculation which is left to
the reader shows that $L \subseteq L_{\beta} + B_{\beta}$.

Now induction hypothesis implies that there exist $n \in \omega$
and $\delta \leq \alpha_1 \leq \alpha_2 \cdots \leq \alpha_n$
such that
\[ \left< H_{\delta} + L_{\beta} \right>_* \subseteq
\sum\limits_{m \leq n}B_{\alpha_m} + H_{\delta}. \] Another
calculation shows that this implies
\[ \left< H_{\delta} +
L \right>_* \subseteq \sum\limits_{m \leq n}B_{\alpha_m} +
B_{\beta} + H_{\delta}. \] This finishes the proof. \qed

{\it (Continuation of the proof of Theorem \ref{main2})} By Lemma
\ref{contained} we can choose $n^* \in \omega$ and $\delta \leq
\alpha_1 \leq \alpha_2 \leq \cdots \leq \alpha_{n^*}$ such that
\[ \left< H_{\delta} + \Z y_{\delta} \right>_* \subseteq
\sum\limits_{m \leq n^*}B_{\alpha_m} + H_{\delta}.
 \]
For every $m \leq n^*$ we choose a finite set $W_m \subset \kappa$
and an integer $n_m \in \omega$ such that
\[ B_{\alpha_m} \subseteq \left< \sum\limits_{\gamma \in W_m}\Z
y_{\gamma} + \sum\limits_{i \leq n_m}x_i \right>_* .\] Collecting
all these generators and letting $W=\bigcup\limits_{m \leq
n^*}W_m$ and $k=\max\{ n_m : m\leq n^* \}$ we obtain
\begin{equation}
\label{nonB2}
\left< H_{\delta} + \Z y_{\delta} \right>_*
\subseteq B + H_{\delta} \end{equation} where $B=\left<
\sum\limits_{\gamma \in W}\Z y_{\gamma} + \sum\limits_{i \leq
k}x_i \right>_*$.

Now choose $\beta \in \delta \backslash W$ and let $n \geq k$
such that
\[ n \in A_{\beta} \cap A_{\delta} \backslash
\bigcup\limits_{\gamma \in W, \gamma \not=\delta } A_{\gamma}. \]
Note that this choice is possible by Lemma \ref{reals}. It is now
straightforward to see that $p_n^{-1}(y_{\delta}-y_{\beta})$ is
an element of $\left< H_{\delta} + \Z y_{\delta} \right>_*$ but
it is not an element of $B + H_{\delta}$ contradicting equation
(\ref{nonB2}). This finishes the proof of Theorem \ref{main2} and
therefore the proof of our Main Theorem \ref{main}. \qed


\goodbreak


\begin{thebibliography}{99}
\bibitem{AH}{\bf U. Albrecht and P. Hill}, {\it Butler groups of infinite rank and axiom 3},
Czech. Math. J., {\bf 37}, (1987), 293--309.

\bibitem{AV}{\bf D.M. Arnold and C. Vinsonhaler}, {\it Finite rank
Butler groups: a survey of recent results}, Lecture Notes in Pure
Appl. Math., {\bf 146} (Marcel Dekker), (1993), 17--41.

\bibitem{Bi}{\bf L. Bican}, {\it Infinite rank Butler groups}, Advances in algebra and model theory, Gordon and Breach Publishers, Algebra Log. Appl. {\bf 9} (1997), 287--317.


\bibitem{B}{\bf M.C.R. Butler}, {\it A class of torsion-free abelian groups of finite rank}, Proc. London Math, Soc. {\bf 15} (1965), 680--698.

\bibitem{BF}{\bf L. Bican and L. Fuchs}, {\it Subgroups of Butler
groups}, Comm. in Algebra, {\bf 22 }, (1994), 1037--1047.

\bibitem{BS}{\bf L. Bican and L. Salce}, {\it Butler groups of infinite rank}, Abelian Group Theory, Lecture Notes in Math. {\bf 1006}, Springer Verlag (1965), 680--698.

\bibitem{DHR}{\bf M. Dugas, P. Hill and K.M. Rangaswamy}, {\it Infinite rank Butler groups
II}, Trans. Amer. Math. Soc., {\bf 320}, (1990), 643--664.

\bibitem{DT}{\bf M. Dugas and B. Thom\'{e}}, {\it The functor Bext under the negation of CH}, Forum Math. {\bf 3} (1991), 23--33.

\bibitem{EM}{\bf P.C. Eklof and A.H. Mekler}, {\it Almost Free
Modules - Set-Theoretic Methods}, North Holland Mathematical
Library, {\bf 46}, (1990).

\bibitem{ER}{\bf P. Erd\"os and R. Rado}, {\it A partition calculus in set theory}, Bull. Am. Math. Aoc., {\bf
62}, (1956), 427--489.

\bibitem{Fu}{\bf L. Fuchs}, {\it Infinite Abelian Groups}, Vol. I and II,
Academic Press (1970 and 1973).


\bibitem{F1}{\bf L. Fuchs}, {\it A survey on Butler groups of infinite rank}, Contemp. Math. {\bf 171} (1994), 121--139).

\bibitem{F2}{\bf L. Fuchs}, {\it Butler groups of infinite rank},
J. Pure Appl. Algebra, {\bf 98}, (1995), 25--44.


\bibitem{FM}{\bf L. Fuchs and M. Magidor}, {\it Butler groups of arbitrary cardinality}, Israel J. Math. {\bf 84} (1993), 239--263.

\bibitem{FR}{\bf L. Fuchs and K.M. Rangaswamy}, {\it Butler groups that are unions of subgroups with countable typesets}, Arch. Math. {\bf 61} (1993), 105--110.


\bibitem{H}{\bf R. Hunter}, {\it Balanced subgroups of abelian groups}, Trans. of the American Math. Soc. {\bf 215} (1976), 81--98.


\bibitem{J}{\bf T. Jech}, {\it Set Theory}, Academic Press, New York (1973).

\bibitem{K}{\bf K. Kunen}, {\it Set Theory - An Introduction to Independent Proofs}, Studies in Logic and the Foundations of Mathematics, North Holland, {\bf 102} (1980).

\bibitem{MS}{\bf M. Magidor and S. Shelah}, {\it $\Bext^2(G,T)$ can be non trivial even assuming GCH}, Contemp. Math. {\bf 171} (1994), 287--294.

\bibitem{R}{\bf K.M. Rangaswamy}, {\it A homological characterization of abelian $B_2$-groups}, Proc. Amer. Math. Soc. {\bf 121} (1994), 409--415.

\bibitem{S}{\bf S. Shelah}, {\it Proper and Improper Forcing}, Perspectives in Mathematical Logic, Springer Verlag (1998).

\end{thebibliography}
\end{document}